\documentclass[aap,MSNbibl,seceqn,dvips]{arximspdf}

\doi{10.1214/12-AAP891} 
\volume{23}
\issue{5}
\pubyear{2013}
\firstpage{1967}
\lastpage{1987}

\makeatletter

\newtheorem{Theorem}{Theorem}[section] 
\newtheorem{MainTheorem}{Theorem}

\newtheorem{MainProposition}{Proposition} 

\newtheorem{MainCorollary}{Corollary}
\newtheorem{claim}[Theorem]{Claim}

\newproclaim{remark}[Theorem]{Remark}
\newproclaim{definition}[Theorem]{Definition}

\newcommand{\cD}{\mathcal D}
\newcommand{\cF}{\mathcal F}
\newcommand{\cG}{\mathcal G}
\newcommand{\cK}{\mathcal K}
\newcommand{\cL}{\mathcal L}
\newcommand{\cN}{\mathcal N}

\newcommand{\bbN}{{\mathbb N}}
\newcommand{\bbP}{{\mathbb P}}
\newcommand{\bbR}{{\mathbb R}}
\newcommand{\bbT}{{\mathbb T}}
\newcommand{\bbZ}{{\mathbb Z}}

\renewcommand{\b}{\beta}
\renewcommand{\c}{\chi}
\newcommand{\g}{\gamma}
\newcommand{\h}{\eta}
\newcommand{\s}{\sigma}
\newcommand{\D}{\Delta}
\renewcommand{\L}{\Lambda}
\renewcommand{\O}{\Omega}

\newcommand{\neper}{e}
\newcommand{\sset}{\subset}
\def\tc{\mid}

\def\Var{\operatorname{Var}}
\def\gap{\operatorname{gap}}

\makeatother

\begin{document}
\begin{frontmatter}

\title{Kinetically constrained spin models on trees\thanksref{T1}}
\runtitle{Kinetically constrained spin models on trees}

\thankstext{T1}{Supported by the European Research Council through
the ``Advanced Grant'' PTRELSS 228032.}

\begin{aug}
\author[A]{\fnms{F.} \snm{Martinelli}\corref{}\ead[label=e1]{martin@mat.uniroma3.it}}
\and
\author[B]{\fnms{C.} \snm{Toninelli}\thanksref{t2}\ead[label=e2]{cristina.toninelli@upmc.fr}}
\runauthor{F. Martinelli and C. Toninelli}
\affiliation{University of Roma Tre and CNRS, University Paris VI--VII}
\address[A]{Dipartimento Matematica\\
Universit\`{a} Roma Tre\\
Largo S.L. Murialdo 00146, Roma\\
Italy\\
\printead{e1}}
\address[B]{Laboratoire de Probabilit\'es\\
\quad et Mod\`eles Al\`eatoires\\
CNRS-UMR 7599 Universit\'es Paris VI-VII 4\\
Place Jussieu F-75252 Paris Cedex 05\\
France\\
\printead{e2}} 
\end{aug}

\thankstext{t2}{Supported in part by the French Ministry of Education
through the ANR-2010-BLAN-0108.}

\received{\smonth{2} \syear{2012}}
\revised{\smonth{8} \syear{2012}}

%
\begin{abstract}
We analyze kinetically constrained $0$--$1$ spin models (KCSM) on rooted
and unrooted trees of finite connectivity.
We focus in particular on the class of Friedrickson--Andersen
models FA-jf and on an oriented version of them.
These tree models
are particularly relevant in physics literature since some of them
undergo an ergodicity breaking transition with the mixed first-second
order character of the glass transition.
Here we first identify the ergodicity regime and prove that the
critical density for
FA-jf and OFA-jf models coincide with that of a suitable bootstrap
percolation model. Next we prove for the first time positivity of the
spectral gap in the whole ergodic regime via a
novel argument based on martingales ideas. Finally, we discuss how
this new technique can be generalized to analyze KCSM on the regular
lattice $\bbZ^d$.
\end{abstract}

%

\end{frontmatter}

\section{Introduction}

Facilitated or kinetically constrained spin models\break (KCSM) are
interacting particle systems which have been introduced in physics
literature~\cite{Fredrickson1,Fredrickson2} to model liquid/glass
transition and more
generally ``glassy dynamics''~\cite{Ritort-Sollich,GST}. They are defined
on a locally finite, bounded degree, connected graph $\mathcal
G=(V,E)$ with vertex set $V$ and edge set $E$. Here we will focus
on models for which the graph is an infinite, rooted or unrooted
tree of finite connectivity $k+1$, which we will
denote by $\bar\bbT^k$ and $\bbT^k$, respectively.
A configuration is given by
assigning to each site $x\in V$ its occupation variable
$\eta_x\in\{0,1\}$ which corresponds to an empty or filled site,
respectively. The evolution is given by a Markovian stochastic dynamics
of Glauber type. Each site\vadjust{\goodbreak} waits an independent, mean one, exponential
time and then, provided the current configuration around it satisfies an
a priori specified constraint, its occupation variable is refreshed to
an occupied or to an empty state with probability $p$ or $1-p$,
respectively. For each site $x$ the corresponding constraint does not
involve $\h_x$, thus detailed balance w.r.t. Bernoulli($p$) product
measure $\mu$ can be easily verified and the latter is an invariant
reversible measure for the process.

Among the most studied KCSM we recall FA-jf models~\cite{Fredrickson1}
for which the constraint requires at least $j$ (which is sometimes
called ``facilitating parameter'') empty sites among the nearest
neighbors.
FA-jf models display a feature which is common to all KCSM introduced
in physics literature: for each vertex $x$ the constraint
imposes a maximal number of occupied sites in a proper neighborhood of
$x$ in order to allow the moves.
As a consequence the dynamics becomes slower at higher density and an
ergodicity breaking transition may occur at a finite critical density
$p_c<1$. This threshold corresponds to the lowest density at which a
site belongs with positive probability to an infinite cluster of
particles which
are mutually and forever blocked due to the constraints; see Section
\ref{ergodicity}.

The FA-jf models on $\mathbb Z^d$ do not display an ergodicity
breaking transition at a nontrivial critical density, that is, $p_c=0$
for $j>d$ and $p_c=1$ otherwise~\cite{CMRT}. On the other hand they do
display such a transition
on nonrooted trees when $1<j<k$~\cite{chalupa,sellitto,schwartz}.
Furthermore if $j\neq k-1$, this transition is expected to display a
mixed first/second character and to share similar features to the mode
coupling transition, a property which makes them particularly
interesting from the point of view of the glass transition~\cite{sellitto}.

Another key feature of
KCSM is the existence of blocked configurations, namely
configurations with all creation/destruction rates identically equal to
zero. This implies the existence of several invariant measures and the
occurrence of unusually long mixing times compared to high-temperature
Ising models (see Section~7.1 of~\cite{CMRT}). Furthermore the
constrained dynamics is usually not attractive so that monotonicity
arguments valid, for example, ferromagnetic stochastic Ising models
cannot be
applied.

Due to the above properties the basic issues concerning the large time
behavior of the process are nontrivial. The first rigorous results
were derived in~\cite{Aldous} for the East model which is defined on
$\mathbb Z$ with the constraint requiring the nearest neighbor site to
the right to be empty. In~\cite{Aldous} it was proven that the
spectral gap of East is positive for all $p<1$ and also that it
shrinks faster than any polynomial in $(1-p)$ as $p\uparrow1$. In
\cite{CMRT} positivity of the spectral gap of KCSM inside the
ergodicity region (i.e., for $p<p_c$) has been proved
in much greater generality and (sometimes sharp) bounds for $p\nearrow
p_c$ were established. These results include FA-jf models on any
$\mathbb Z^d$ for any choice of the facilitating parameter $j$ and of
the spatial dimension~$d$.

The technique developed in~\cite{CMRT} cannot be applied to models
on trees because of the exponential growth of the
number of vertices and, so far, very few rigorous results have been
established.
Indeed the only models for which results on the spectral gap are
available are:
(i) the FA-1f model on $\bbT^k$ and $\bar\bbT^k$ (actually on a
generic connected graph) and (ii) the so-called East model on
$\bar\bbT^k$ for which the root is unconstrained while, for any other
vertex $x$, the constraint requires the ancestor of $x$ to be empty.
For these specific models $p_c=1$ and the positivity of the spectral
gap has been proven in~\cite{LNM} in the whole ergodicity region and
for any choice of the graph connectivity.

Here we will study FA-jf models on $\bbT^k$ and $\bar\bbT^k$ for
$1<j\leq k$
together with
a new class of models that we call oriented FA-jf models (OFA-jf). In the
OFA-jf model the constraint at $x$ requires at least $j$ empty sites
among the \textit{children} of $x$.

We first prove that the ergodicity threshold $p_c$ for the FA-jf and
OFA-jf models, with the same choice for the parameter $j$ and the same
graph connectivity $k+1$, coincide and it is nontrivial (see Theorem
\ref{teopc}).
Then we prove positivity of the spectral gap in the whole ergodicity regime
for the oriented OFA-jf models. Finally, by combining the above
results together with an appropriate comparison technique, we
establish positivity of the spectral gap in the whole ergodicity
regime for the FA-jf models. The results concerning the spectral gap
can be found in Theorem~\ref{teogap} and a simple application to
the mixing time of finite system in Corollary~\ref{Coro}. Finally, in
the nonergodic
regime, we prove that, for the oriented or nonoriented FA-jf models,
the spectral gap shrinks to zero exponentially fast in the system
size; see Theorem~\ref{abovepc}.

The new technique devised to study constrained models on trees can be
generalized to deal also with KCSM on other graphs.
In Section~\ref{generalizations} we discuss how one can recover the
result of positivity of the spectral gap in the ergodic regime for
models on
$\bbZ^d$. We detail in particular the case of the north--east model on
$\bbZ^2$ (Theorem~\ref{gapNE}), a result which was already derived in
\cite{CMRT} but with a completely different (and more lengthy) technique.

\section{Models and main results}
\label{main}
\subsection{Setting and notation}

\subsubsection*{The graphs}
The models we consider are either defined on the infinite regular tree
of connectivity $k+1$, in the sequel denoted by $\bbT^k$ or on the
infinite, rooted
$k$-ary tree $\bar\bbT^k$.
In the unrooted case each vertex $x$ has $k+1$ neighbors, while in
the rooted case each vertex different from the root has $k$ children
and one ancestor, and the root $r$ has only $k$ children. In the sequel
we will denote by $V$ the set of vertices of either $\bbT^k$ or
of $\bar\bbT^k$ whenever no confusion arises, by $\mathcal N_x$ the
set of neighbors of a given
vertex $x$ and, in the rooted case, by $\mathcal K_x$ the set of its
children. In the rooted case we denote by $d_x$ the \textit{depth} of
the vertex $x$, that is, the graph distance between $x$ and the root
$r$.

\subsubsection*{The configuration spaces} For both oriented and nonoriented
models we choose as configuration space the set $\Omega=\{0,1\}^V$ whose
elements will usually be assigned Greek letters. We will often write
$\h_x$ for the value
at $x$ of the element $\h\in\O$. We will also write
$\O_A$ for the set $\{0,1\}^A$, $A\subseteq V$. With a slight abuse of
notation, for
any $A\subseteq V$ and any $\eta,\omega\in\O$, we let $\eta_A$ to be
the restriction of $\eta$ to the set $A$ and $\eta_A\cdot\omega_{A^c}$ to be the
configuration which equals $\eta$ on $A$ and $\omega$ on $V\setminus
A$.

\subsubsection*{Probability measures}
For any $A\subseteq V$ we denote
by $\mu_A$ the product measure $\bigotimes_{x\in A} \mu_x$ where each
factor $\mu_x$ is the Bernoulli measure on $\{0,1\}$ with $\mu_x(1)=p$
and $\mu_x(0)=q$ with $q=1-p$. If $A=V$ we abbreviate $\mu_V$ to
$\mu$.

\subsubsection*{Conditional expectations and conditional variances}
Given a function $f\dvtx \Omega\to\bbR$ depending on finitely many variables,
in the sequel referred to as \textit{local function}, and a set $A\sset
V$ we define the
function $\eta\mapsto\mu_A(f)(\eta)$ by the formula
\[
\mu_A(f) (\eta):= \sum_{\sigma\in\Omega_A}
\mu_A(\sigma)f(\sigma_A\cdot\eta_{A^c}).
\]
Clearly $\mu_A(f)$ coincides with the
\textit{conditional expectation} of $f$ given the configuration outside $A$.
Similarly we write $\Var_A(f)=\mu_A(f^2)-\mu_A(f)^2$ for
the \textit{conditional variance} of $f$ given $\eta_{A^c}$. Note that
$\Var_A(f)=0$
if and only if
$f$ does not depend on the configuration inside $A$. In case $A=V$ we
abbreviate $\Var_V(f)$ to $\Var(f)$.

\subsection{Facilitated models}
%
\begin{definition}\label{defFA}
Fix $k\in\bbZ_+$ and a facilitating
parameter $j\in[1,\ldots, k]$. The FA-jf and OFA-jf models at density
$p$ are continuous time Glauber-type Markov
processes\vspace*{1pt} on $\O$, reversible w.r.t. $\mu$,
with Markov semigroups $P_t= \neper^{t\cL}$ and $\bar P_t=
\neper^{t\bar\cL}$,
respectively, whose infinitesimal generators $\cL,\bar\cL$ act
on local functions $f\dvtx \O\mapsto\mathbb R$ as follows:
%
\begin{eqnarray}
\cL f(\omega)&=&\sum_{x\in\bbT^k}c_{x}(\omega)
\bigl[\mu_x(f) (\omega)-f(\omega) \bigr],
\\
\label{FAjf}
\bar\cL f(\omega)&=&\sum_{x\in\bar\bbT^k}\bar c_{x}(
\omega) \bigl[\mu_x(f) (\omega)-f(\omega) \bigr].
\end{eqnarray}
The function $c_x$ (or $\bar c_x$), in the sequel referred to as the
\textit{constraint at $x$}, is defined by
%
\begin{eqnarray}
\label{cx} c_{x}(\omega)&=& %
\cases{1, &\quad if $ \displaystyle \sum
_{y\in\mathcal N_x}(1-\omega_y)\geq j$,
\vspace*{2pt}\cr
0, &\quad otherwise, }
\\
\label{therates}
\bar c_{x}(\omega)&=& %
\cases{ 1, &\quad if $\displaystyle  \sum
_{y\in\mathcal K_x}(1-\omega_y)\geq j$,
\vspace*{2pt}\cr
0, &\quad otherwise. }
\end{eqnarray}
\end{definition}
It is easy to check by standard methods (see, e.g.,~\cite{Liggett}) that
the processes are well defined and that their generators can be
extended to
nonpositive self-adjoint operators on $L^2(\bbT^k,\mu)$ and
$L^2(\bar\bbT^k,\mu)$,
respectively.

Both processes can of course be defined also on finite regular trees,
rooted or unrooted. In this case and in order to ensure irreducibility
of the Markov chain the constraints must be suitably modified.
%
\begin{definition}\label{finite}
Let $\bbT$ be a finite subtree of either $\bbT^k$ or of $\bar\bbT^k$ and let, for any $\eta\in\O_\bbT$, $\eta^0\in\O$ denote the
extension of $\eta$ in $\O$ given
by
\[
\eta^0_x= %
\cases{ \eta_x, &\quad if
$x\in\bbT$,
\cr
0, &\quad if $x\in\bbT^k\setminus\bbT$. } %
\]
For any $x\in\bbT$ define the \textit{finite constraints} $c_{\bbT,x},\bar
c_{\bbT,x}$ by
\[
c_{\bbT,x}(\eta)= c_x\bigl(\eta^0\bigr),\qquad \bar
c_{\bbT,x}(\eta)= \bar c_x\bigl(\eta^0\bigr)
\qquad\forall\eta\in\O_\bbT.
\]
We will then refer to the \textit{OFA-jk
model} or the \textit{FA-jk model} on $\bbT$ as the irreducible,
continuous time
Markov chains on $\O_\bbT$ with generators
%
\begin{eqnarray}
\label{fin-gen2} \cL_\bbT f(\eta)&=&\sum_{x\in\bbT}c_{\bbT,x}
\bigl[\mu_x(f)-f\bigr] \eta \in\O_\bbT,
\\
\bar\cL_\bbT f(\eta)&=&\sum_{x\in\bbT}\bar
c_{\bbT,x}\bigl[\mu_x(f)-f\bigr] \eta\in\O_\bbT,
\end{eqnarray}
respectively.
\end{definition}

\subsection{Ergodicity}
\label{ergo}
Given $k,j\in\bbZ_+$ with $j\le k$, it is natural to define (see~\cite{CMRT})
a critical density for
each model as follows:
%
\begin{eqnarray}
\label{eq1} p_c &=&\sup\bigl\{p\in[0,1]\dvtx \mbox{0 is simple eigenvalue
of } \cL\bigr\},
\\
\bar p_c &=&\sup\bigl\{p\in[0,1]\dvtx \mbox{0 is simple eigenvalue of }
\bar\cL\bigr\}.
\end{eqnarray}
The regime $p<p_c$ or $p<\bar p_c$ is called the {\sl ergodic region}
and we say that an {\sl ergodicity breaking transition} occurs at the
critical density. We will first establish the coincidence of the
critical threshold for oriented and unoriented models.
%
\begin{MainTheorem}
\label{teopc}
Given $k,j\in\bbZ_+$ with $j\le k$, let $g_p(\lambda):=p\sum_{i=k-j+1}^{k}{k\choose i}\lambda^i (1-\lambda)^{k-i}$ and
define
%
\begin{equation}
\label{eq2} \tilde p:=\sup\bigl\{p\in[0,1]\dvtx  \lambda=0 \mbox{ is the unique
fixed point of } g_p(\lambda)\bigr\}.
\end{equation}
Then $p_c=\bar p_c=\tilde p$ and for any $p<\tilde p$
the value $0$ is a simple
eigenvalue of the generators $\cL$ and $\bar\cL$. Moreover $\tilde
p\in(0,1)$ if and only if $2\le j\le k$.
\end{MainTheorem}
We then turn to the study of the relaxation to equilibrium in
$L^2(\mu)$. A~key object here is the spectral gap (or inverse of the
relaxation time) of the generator $\cL$ (or~$\bar\cL$), defined as
%
\begin{equation}
\label{eqgap} \gap(\cL):= \mathop{\inf_{f\in\operatorname{Dom}(\cL)}}_{ f\neq
\mathrm{const}}
\frac{\cD(f)}{\Var(f)},
\end{equation}
where the Dirichlet form $\cD(f)$ is the quadratic form
$\cD(f)=\mu(f,-\cL f)$ associated to $-\cL$.
Indeed a positive spectral gap implies that the reversible measure $\mu
$ is mixing for the
semigroup $P_t$ with exponentially decaying correlations,
\[
\biggl(\int d\mu(\eta) \bigl[P_t f(\eta)-\mu(f)\bigr]^2
\biggr)^{1/2}\leq e^{-\gap(\cL) t}\Var(f) \qquad\forall f\in L^2(
\mu).
\]

\subsection{Main results on relaxation to equilibrium}
For the reader's convenience we split the presentation of our results
into three
sub-sections according to whether $p$ is below, above or equal to the
critical value $p_c$.

\subsubsection{The sub-critical case $p<p_c$}

\begin{MainTheorem}
\label{teogap}
Given $k,j\in\bbZ_+$ with $j\le k$, fix $p<p_c=\bar p_c$. Then
$\gap(\cL)>0$ and $\gap(\bar\cL)>0$.
\end{MainTheorem}
%
\begin{remark}
\label{finite-infinite}
Exactly as in~\cite{CMRT} (see Proposition 2.13 there), in order
to prove positivity of the spectral gap for the infinite trees $\bbT^k$
or $\bar\bbT^k$, it is enough to prove a lower bound on the spectral
gap of the corresponding models on finite balls which is uniform in the
size of the ball.
\end{remark}
It is important to observe that in the oriented case the above
result completes the proof of the exponential decay to equilibrium when
$p<p_c$ and the initial distribution is either a Bernoulli product
measure with density $p'\neq p, p'<p_c$, or it is a $\delta$-measure
on a
deterministic configuration which does not contain blocked
clusters. These results were indeed proven in~\cite{CMST} (see
Theorems 4.2 and 4.3) modulo the hypothesis of positivity of the
spectral gap in the ergodic region.

We finally observe that the above result says nothing about the
behavior of the
spectral gap as a function of $p_c-p$ when $p\uparrow p_c$. See,
however, Section~\ref{critical} below for some work in progress in
this direction.

Our second result, a natural corollary of the spectral gap bounds of
Theorem~\ref{main}, concerns
mixing times of the oriented model on finite sub-trees of $\bar\bbT^k$. In
order to state it we need few extra notation.

Let $\bbT$ be the finite rooted tree consisting of the first $n$
levels of $\bar\bbT^k$.
For any $\eta\in\O_\bbT$ we denote by $\nu_t^\eta$ the law at
time $t$ of the Markov chain with
generator $\bar\cL_{\bbT}$ and by $h_t^\eta$ the relative density
w.r.t. $\mu_\bbT$
of $\nu_t^\eta$, namely
\[
h_t^\eta(\s):=\nu_t^\eta(\sigma)/
\mu_\bbT(\s).
\]
Following~\cite{Gine}, we define the family of mixing times $\{T_a\}_{a\ge1}$ by
\[
T_a:= \inf \Bigl\{t\ge0\dvtx  \max_\eta\mu_\bbT
\bigl(\bigl|h_t^\eta -1\bigr|^a \bigr)^{1/a}
\le1/4 \Bigr\}.
\]
Notice that $T_1$ coincides with the usual mixing time $T_{\mathrm{mix}}$ of
the chain (see,
e.g.,~\cite{Peres}) and that, for any $a\ge1$, $T_1\le T_a$.
%
\begin{MainCorollary}
\label{Coro}
Given $k,j\in\bbZ_+$ with $j\le k$, fix $p<p_c$. Then there
exists a constant $c$ such that
\[
c^{-1}n \le T_1\le T_2 \le c n.
\]
\end{MainCorollary}
%
\begin{remark}
A key ingredient
for the proof of the above Corollary will be the fact that the
marginal of the law $\nu_t^\eta$ over $\O_{\bbT\setminus r}$ is given
by the product of the marginals over the individual subtrees rooted
at the children of the root. Such a property is no longer true in the
unoriented case. In this more complicate setting a possible route to
get a (poorer) bound on the mixing time is the following.

Use a
comparison between the Dirichlet forms of the FA-jk and OFA-jk models to
get that the logarithmic Sobolev constant (see, e.g.,~\cite{Gine}) of
the FA-jk model on a finite regular tree $\bbT\subset\bbT^k$, with
$n$ levels and centered
at a vertex $r$, is bounded from below by constant${}\times{}$the
logarithmic Sobolev constant of the OFA-jk model on the
finite trees $\bbT\setminus r$. Then use the left part of the
well-known bound (see Corollary~2.2.7 in
\cite{Gine})
\begin{eqnarray*}
&&
(\mbox{log-Sobolev constant})^{-1} \\
&&\qquad\le T_2 \le
\mathrm{const}\times (\mbox{log-Sobolev constant})^{-1}\log \bigl(\bigl|\log
\bigl(\mu_\bbT^*\bigr)\bigr| \bigr),
\end{eqnarray*}
where $\mu_\bbT^*:=\min_\eta\mu_\bbT(\eta)$ to infer that the
logarithmic Sobolev constant of the OFA-jk model is
bounded from below by $\mathrm{const}\times T^{-1}_2$. Hence the
logarithmic Sobolev constants of
both the OFA-jk and the FA-jk models on $\bbT$ are bounded from below
by $\mathrm{const} \times
n^{-1}$. Finally use the right part of the above bound to
conclude that the
mixing time $T_2$ for the FA-jk model on $\bbT$ is $O(n^2)$.
\end{remark}

\subsubsection{The super-critical phase $p>p_c$}

Our first result roughly says that, when $p>p_c$, the occupation
number for
the process defined on the infinite
tree does not equilibrate in $L^2(\mu)$.

Denote by $r$ either the
root (in the oriented case) or an arbitrary vertex of $\bbT^k$ (in the
unoriented case).
%
\begin{MainProposition}
\label{soprapc}
Given $k,j\in\bbN$ with $j\le k$, fix $p>p_c$. Then
\[
\lim_{t\to\infty}\Var(\bar P_t \eta_r)>0,
\]
and the same inequality holds with $P_t$ instead of $\bar P_t$.
\end{MainProposition}
The second result concerns the spectral gap on finite balls. Given
$n\in\bbZ_+$ and $r\in\bbT^k$, denote by $\bbT$ either the ball in
$\bbT^k$ of
radius $n$ and center $r$ (in the unoriented case) or the rooted tree
consisting of the first
$n$ levels of $\bar\bbT^k$ (in the oriented case).
%
\begin{MainTheorem}
\label{abovepc}
Given $k,j\in\bbN$ with $j\le k$, fix $p>p_c$. Then there exists
$c>0$ such that
\begin{eqnarray*}
e^{-cn} &\le&\gap(\cL_\bbT) \le e^{-n/c},
\\
e^{-c n} &\le&\gap(\bar\cL_\bbT) \le e^{-n/c}.
\end{eqnarray*}
\end{MainTheorem}
%
\subsubsection{The critical phase $p=p_c$}
\label{critical}

The critical case $p=p_c$ is much more delicate and a detailed analysis
is postponed
to future work~\cite{critical-tree}. We anticipate here that it is possible
to show that the spectral gap on a ball of radius $n$ shrinks at least
polynomially fast in $n^{-1}$. In the rooted case with $j=k$ one can
also prove a converse poly($1/n$) \textit{lower bound} (a much harder
task). These two results then imply that in the rooted case and for $j=k$,
there
exist three positive constants $\beta\ge2$, $c_1,c_2$ such that, for
$p_c-p\ll1$,
\[
c_1 (p_c-p)^\b\le\gap(\bar\cL)\le
c_2(p_c-p)^{2}.
\]
If $2\le j<k$, the analysis of the lower bound on the spectral gap
becomes much more difficult because of
the \textit{discontinuous} character of the bootstrap percolation
transition. More precisely, and contrary to what happens for $j=k$, for
$p=p_c$ the root $r$ belongs to an infinite
blocked cluster with \textit{positive} probability. In this
case it is still unclear whether a poly($1/n$)
lower bound on the spectral gap still holds.

\section{\texorpdfstring{Ergodicity threshold and blocked clusters: Proof of Theorem \protect\ref{teopc}}
{Ergodicity threshold and blocked clusters: Proof of Theorem 1}}
\label{ergodicity}
%
\begin{definition}
Given $k,j\in\bbZ_+$ with $j\le k$, the bootstrap map $B\dvtx\break \{0,1\}^{\bbT^k}\to
\{0,1\}^{\bbT^k}$ associated to the FA-jf model is defined by
%
\begin{equation}
\label{eqbootstrapmap} B(\h)_x=0 \quad\mbox{if either}\quad
\h_x=0 \quad\mbox{or}\quad c_x(\h)=1
\end{equation}
with $c_x$ defined in (\ref{cx}). Analogously we define the bootstrap
map $\bar B$ for the OFA-jf model by replacing $c_x$ with $\bar c_x$ of
(\ref{cx}).
\end{definition}
Having defined the bootstrap map $B$ it is natural to denote by $\mu^{(n)}$
the probability
measure obtained by iterating $n$-times the map $B$\vadjust{\goodbreak}
starting from~$\mu$. In other words, for any $A\sset\O$ $\mu^{(n)}(A)= \mu (\eta\dvtx
B^n(\eta)\in A )$.
As $n$ tends to infinity $\mu^{(n)}$ converge to
a limiting measure $\mu^{(\infty)}$~\cite{Schonmann}, and it is natural
to define the bootstrap percolation threshold $p_{\mathrm{bp}}$ as
the supremum
of the density $p$ of $\mu$ such that $\mu^{(\infty)}$ is concentrated
on the empty configuration.
Analogously we can define $\bar\mu^{(n)}, \bar\mu^{(\infty)}$ and
$\bar p_{\mathrm{bp}}$ in the oriented case.

It is quite clear that the two thresholds $p_{\mathrm{bp}}$ and $\bar p_{\mathrm{bp}}$
must coincide. Choose in fact an arbitrary vertex $r\in\bbT^k$ and write
the unrooted tree $\bbT^k$ as $\bbT^k=\{r\}\bigcup_{y\in\cN_r}\bar
\bbT^k_y$ where each $\bar\bbT^k_y$ is a copy of $\bar\bbT^k$ with
root at $y$. If $p<\bar p_{\mathrm{bp}}$, then a.s. each $y\in\cN_r$ becomes
eventually empty under the bootstrap map $\bar B$ applied to $\bar
\bbT_y^k$ and therefore also under
the less-restrictive bootstrap map~$B$. Thus $p\le p_{\mathrm{bp}}$. On the other
hand, when $p>\bar p_{\mathrm{bp}}$ the set
\[
\cG=\bigl\{\eta\in\O\dvtx  \eta_r=1 \mbox{ and } (\bar
B)^\infty(\eta_{\bar\bbT_y^k})_r=1\ \forall y\in
\cN_r\bigr\}
\]
has positive probability and moreover $B^\infty(\eta)_r=1$ for any
$\eta\in
\cG$. Hence $p\ge p_{\mathrm{bp}}$.

That $p_{\mathrm{bp}}$ coincide
with the third threshold $\tilde p$ given in (\ref{eq2}) has been
established in Proposition 1.2 of~\cite{Peres} (see also
\cite{sellitto,FS} and~\cite{STBT} for an
extension to hyperbolic lattices).
For completeness we shortly reprove this result by showing
that $\bar p_{\mathrm{bp}}=\tilde p$.

We first observe that $\bar\mu^{(\infty)}(\h_r=1)=0$ if and only if
$\lim_{n\to
\infty}\bar p_n=0$ where $\bar p_n:=\bar
\mu^{(n)} (\h_r=1 )$. Second one easily checks that the
nonincreasing
sequence
$\{\bar p_n\}_{n\ge0}$
obeys the recursive equation
$\bar p_n = g_p(\bar p_{n-1})$
with initial condition \mbox{$\bar p_0=p$}. Here $g_p(\cdot)$ has the expression
\[
g_p(\lambda):=p\sum_{i=k-j+1}^{k}
\pmatrix{k
\cr
i}\lambda^i (1-\lambda)^{k-i}.
\]
We now claim that $\lim_{n\to
\infty}\bar p_n=0$ if and only if $p<\tilde p$. In
order to prove the claim we first observe that $\lim_{n\to
\infty}\bar p_n $ is a fixed point of the map $g_p$ and that it is
a nondecreasing function of $p$. Hence $p<\tilde p\Rightarrow\lim_{n\to
\infty}\bar p_n=0$.

To prove the converse we compute
\begin{eqnarray*}
\frac{d}{d\lambda} g_p(\lambda)&=& p\sum_{i=k-j+1}^{k}
\pmatrix{k
\cr
i} \bigl[ i \lambda^{i-1} (1-\lambda)^{k-i}-(k-i)
\lambda^i(1-\lambda)^{k-i-1} \bigr]
\\
&=& p \Biggl[ \sum_{i=k-j}^{k-1}k\pmatrix{k-1
\cr
i} \lambda^{i} (1-\lambda)^{k-1-i}\\
&&\hspace*{10.5pt}{} - \sum
_{i=k-j+1}^{k-1}k\pmatrix{k-1
\cr
i} \lambda^{i}
(1-\lambda)^{k-1-i} \Biggr]
\\
&=& pk \bbP(N_{\lambda,k}=k-j)>0,
\end{eqnarray*}
where $N_{\lambda,k}\sim\operatorname{Binom}(k-1,\lambda)$.\vadjust{\goodbreak}

Therefore $g_p$ is strictly increasing in $(0,1)$, and if it has a
fixed point $\lambda^*\in(0,p)$, then
necessarily $\lim_{n\to
\infty}\bar p_n \ge\lambda^*$. Hence
$\lim_{n\to
\infty}\bar p_n=0\Rightarrow p<\tilde p$.

We finally check that $\tilde p\in(0,1)$ if and only if $2\le j\le k$. The
Markov inequality implies that
\[
g_p(\lambda)\le p\frac{k}{k-j+1}\lambda.
\]
Hence
$g_p(\lambda)<\lambda$ if $j=1$ and $p<1$. When $j\in[2,k]$ it is
also clear that
$\tilde p\in[\tilde p_2,\tilde p_k]$, where $\tilde p_2,\tilde p_k$
correspond to the extreme cases $j=2$ and $j=k$, respectively. When
$j=k$ the threshold $\tilde p_k$ coincides with the usual site
percolation threshold $1/k$ (see~\cite{Grimmett}). When $j=2$ and $k\ge3$ an exact
computation~\cite{Peres} gives
\[
\tilde p_2= \frac{(k-1)^{2k-3}}{k^{k-1}(k-2)^{k-2}} <1.
\]

\begin{remark}
\label{bdecay}
It is not difficult to check that, for $k\ge2$, the limit as $n\to
\infty$ of both
sequences $ \{\mu^{(n)}(\eta_x=1) \}_{n\ge0}$ and $ \{
\bar
\mu^{(n)}(\eta_x=1) \}_{n\ge0}$ is:
\begin{itemize}
\item zero and attained at least exponentially fast if $p<\tilde p$;
\item zero and attained polynomially fast (in
$1/n$) for $j=k$ and $p=\tilde p$;
\item strictly positive for $j\in[2,k)$ and $p=\tilde p$.
\end{itemize}
\end{remark}

The proof of Theorem~\ref{teopc} now follows from the above discussion
together with the following
proposition which can be proved following exactly the same lines as
Proposition 2.5 of~\cite{CMRT}.
%
\begin{MainProposition}\label{propo}
$p_c=p_{\mathrm{bp}}$ and $\bar p_c=\bar p_{\mathrm{bp}}$.
\end{MainProposition}

\section{Relaxation to equilibrium: Proofs}

\subsection{The sub-critical phase $p<p_c$}
\label{subcrit}
In what follows we fix once and for all $j,k\in\bbZ_+$ with $j\le k$,
together with a density $p\in[0,p_c)$.
\begin{pf*}{Proof of Theorem~\ref{teogap}: the oriented case}
We begin by proving positivity of the
spectral gap in the oriented case OFA-jf at density $p$.

We first fix some additional notation. We denote
by $\bbT$ the finite $k$-ary tree consisting of the first $n$ levels
(counting the root $r$)
of $\bar\bbT^k$, where $n$ should be thought of as arbitrarily large
compared to all
other constants. For $x\in\bbT$, $\bbT_x$ will
denote the $k$-ary sub-tree of $\bbT$ rooted at $x$ and with $n-d_x+1$
levels, where $d_x\in[1,n]$ is the level label of $x$. We also set
$\hat\bbT_x:=\bbT_x\setminus\{x\}$. In the sequel we shall refer to
the number of levels $n-d_x+1$ as the \textit{depth} of the tree $\bbT_x$.


The key idea for the proof is to introduce long-range constraints.
%
\begin{definition}
\label{meso}
For any $\eta\in\O_{\hat\bbT_x}$, let $\eta^1\in\O_\bbT$ be
equal to $\eta$
in $\hat\bbT_x$ and equal to $1$ in $\bbT\setminus
\hat\bbT_x $. Then, for any integer $\ell$ we define
\[
\bar c_x^{(\ell)}(\h)= %
\cases{1, &\quad if the depth
of $\bbT_x$ is not larger than $\ell$ or if $ (\bar B
)^\ell\bigl(\eta^1\bigr)_x=0$,
\cr
0, &\quad
otherwise. } %
\]
\end{definition}
%
In what follows we will first consider an auxiliary long-range,
kinetically constrained model on $\bbT$ whose infinitesimal generator
is as in (\ref{fin-gen2})
but with $\bar c_{\bbT,x}$ substituted by
$\bar c_x^{(\ell)}$. We will show that this auxiliary model
has
a spectral gap which is bounded away from zero \textit{uniformly} in the depth
$n$ of $\bbT$, provided $\ell$ is large enough depending on
$p,j,k$. Then we will apply standard
comparison arguments between the Dirichlet forms with constraints
$\bar c_{\bbT,x}$ and $\bar c_x^{(\ell)}$ to show that also the
original model has a
spectral gap which is uniformly positive in $n$. By appealing to Remark
\ref{finite-infinite}
that completes the proof.

Let $\cD_\bbT^{(\ell)}(f)$ denote the new Dirichlet form
corresponding to
the generator
\[
\cL_\bbT^{(\ell)} f(\omega)=\sum_{x\in\bbT}
\bar c^{(\ell
)}_{x}(\omega) \bigl[\mu_x(f)-f(\omega)
\bigr]
\]
with
the auxiliary constraints $\bar c_x^{(\ell)}$, that is,
\[
\cD_\bbT^{(\ell)}(f)= \frac12 \sum
_{x\in\bbT}\mu_\bbT \bigl(\bar c_x^{(\ell)}
\Var_x(f) \bigr).
\]
Our aim is to establish the so-called \textit{Poincar\'e inequality}
%
\begin{equation}
\label{PI} \Var_\bbT(f)\le\lambda\cD_\bbT^{(\ell)}(f)
\qquad\forall f\dvtx \O_\bbT\mapsto\bbR
\end{equation}
for some constant $\lambda$ independent of the depth $n$ of the tree
$\bbT$.
%
\begin{remark}
Notice that (\ref{PI}) is the natural analog of the renormalized
Poincar\'e inequality in~\cite{CMRT}; see formula (5.1) there.
\end{remark}
For the reader's convenience we begin by recalling some elementary
properties of
the variance which will be applied in the sequel.
Consider two probability spaces
$ (\O_i,\cF_i,\nu_i )$, $i=1,2$, together with their product
probability space $ (\O,\cF,\nu )$. Then, for any $f\in
L^2(\O,\nu)$,
\[
\Var(f) \le\nu \bigl(\Var(f\tc\cF_1)+ \Var(f\tc\cF_2)
\bigr) \quad\mbox{and}\quad \Var\bigl(\nu(f\tc\cF_2)\bigr)\le\nu_1
\bigl(\Var(f\tc \cF_1) \bigr)
\]
so that
%
\begin{equation}
\label{eq4} \Var(f)\le\nu \bigl(\Var(f\tc\cF_1)+\Var(f\tc
\cF_2) \bigr).
\end{equation}
Clearly $\Var(f\tc\cF_1)=\nu_2(f^2)-\nu_2(f)^2$, $\nu(f\tc\cF_2)=\nu_1(f)$ and so forth. Moreover,
%
\begin{equation}
\label{var} \Var(f)= \nu \bigl(\Var(f\tc\cF_2) \bigr)+ \Var \bigl(
\nu(f\tc \cF_2) \bigr).
\end{equation}
Back to the proof and motivated by~\cite{MSW} we first claim that
%
\begin{equation}
\label{A} \Var_\bbT(f)\leq\sum_{x\in\bbT}
\mu_\bbT \bigl(\Var_x \bigl(\mu_{\hat\bbT_x}(f) \bigr)
\bigr).
\end{equation}
To prove the claim we proceed recursively on the depth $n$ of
$\bbT$. The claim is trivially true for $n=0$. We now assume (\ref
{A}) when $\bbT$ has depth $n-1$, and using the formula
for the conditional variance we write
%
\begin{equation}
\label{B} \Var_\bbT(f) = \mu_\bbT \bigl(
\Var_\bbT(f\tc\h_r) \bigr) + \Var_\bbT \bigl(
\mu_\bbT(f\tc\h_r) \bigr).
\end{equation}
Notice that, given the spin $\h_r$ at the root, $\Var_\bbT(f\tc\h_r)$ is
nothing but the variance of $f$ w.r.t. the product measure
$\mu_{\bbT\setminus\{r\}}=\prod_{y\in\cK_x}\mu_{\bbT_y}$. Thus
\[
\Var_\bbT(f\tc\h_r)\le\sum_{y\in\cK_x}
\mu_\bbT \bigl(\Var_{\bbT_y}(f)\tc\h_r \bigr)
\]
and
\[
\mu_\bbT \bigl(\Var_\bbT(f\tc\h_r) \bigr)\le
\sum_{y\in\cK
_x}\mu_\bbT \bigl(
\Var_{\bbT_y}(f) \bigr).
\]
Each one of the sub-trees $T_y$ has depth $n-1$, and therefore the
inductive assumption implies that
%
\begin{eqnarray}\label{C}
\sum_{y \in\cK_x}\mu_\bbT \bigl(
\Var_{\bbT_y}(f) \bigr) &\le& \sum_{y \in\cK_x}\sum
_{z\in\bbT_y}\mu_\bbT \bigl(\Var_{z}
\bigl(\mu_{\hat T_z}(f) \bigr) \bigr)
\nonumber\\[-8pt]\\[-8pt]
&=& \mathop{\sum_{x\in\bbT}}_{x\neq r} \mu_\bbT \bigl(
\Var_x \bigl(\mu_{\hat\bbT_x}(f) \bigr) \bigr).\nonumber
\end{eqnarray}
By putting together the right-hand side of (\ref{C}) with the last
term in (\ref{B}), we get the claim for depth $n$.

We now examine a generic term $\mu_\bbT (\Var_x (\mu_{\hat\bbT_x}(f) ) )$ in the right-hand side of (\ref
{A}). We write
\[
\mu_{\hat\bbT_x}(f)=\mu_{\hat\bbT_x}\bigl(\bar c_x^{(\ell)}f
\bigr) + \mu_{\hat\bbT_x}\bigl(\bigl[1-\bar c_x^{(\ell)}
\bigr]f\bigr)
\]
so that
%
\begin{equation}\label{D0}
\Var_x \bigl( \mu_{\hat\bbT_x}(f) \bigr) \le2 \Var_x
\bigl(\mu_{\hat\bbT_x} \bigl(\bar c_x^{(\ell)}f \bigr)
\bigr) + 2 \Var_x \bigl(\mu_{\hat\bbT_x} \bigl(\bigl(1-\bar
c_x^{(\ell)}\bigr)f \bigr) \bigr).
\end{equation}
The Cauchy--Schwarz inequality shows that
%
\begin{equation}\label{D}
\Var_x \bigl(\mu_{\hat\bbT_x} \bigl(\bar c_x^{(\ell)}f
\bigr) \bigr)\le\mu_{\hat\bbT_x} \bigl(\Var_x \bigl(\bar
c_x^{(\ell
)}f \bigr) \bigr)= \mu_{\hat\bbT_x} \bigl(\bar
c_x^{(\ell)}\Var_x(f) \bigr),
\end{equation}
because $\bar c_x^{(\ell)}$ does not depend on the spin at $x$. Notice
that the right-hand side in (\ref{D}) is just the contribution of the
root to the Dirichlet form $\cD_\bbT^{(\ell)}(f)$.

We now turn to the analysis of the more complicated second term\break $\Var_x (\mu_{\hat
\bbT_x} ((1-\bar c_x^{(\ell)})f ) )$, in the
nontrivial case
$n-d_x+1>\ell$. We write
%
\begin{eqnarray}\qquad
\Var_x \bigl(\mu_{\hat\bbT_x} \bigl(\bigl(1-\bar
c_x^{(\ell)}\bigr)f \bigr) \bigr)&=&\Var_x \bigl(
\mu_{\hat\bbT_x} \bigl(\bigl(1-\bar c_x^{(\ell
)}\bigr)
\bigl(f-\mu_{ \bbT_x}(f)+\mu_{ \bbT_x}(f) \bigr) \bigr)\bigr)
\nonumber\\[-8pt]\\[-8pt]
&=&\Var_x \bigl(\mu_{\hat\bbT_x} \bigl(\bigl(1-\bar
c_x^{(\ell
)}\bigr)g \bigr) \bigr),\nonumber
\end{eqnarray}
where $g:=f-\mu_{ \bbT_x}(f)$ and we use the fact that $\mu_{\hat
\bbT_x}((1-\bar c_x^{(\ell)})\mu_{ \bbT_x}(f))$ does not depend on
$\eta_x$.
Recall that the constraint $\bar c_x^{(\ell)}$ depends only on the
spin configuration in the first $\ell$ levels below $x$, in the sequel
denoted by $\D_x$. Then
%
\begin{eqnarray}
\label{least} \Var_x \bigl(\mu_{\hat\bbT_x} \bigl(\bigl(1-\bar
c_x^{(\ell)}\bigr)g \bigr) \bigr)&\leq&\mu_x \bigl(
\bigl(\mu_{\hat\bbT_x} \bigl(\bigl(1-\bar c_x^{(\ell)}\bigr)
\mu_{\hat\bbT_x\setminus\Delta_x }g \bigr) \bigr)^2 \bigr)
\nonumber\\
&\leq&\mu_x \bigl(\mu_{\hat\bbT_x}\bigl(1-\bar
c_x^{(\ell
)}\bigr)\mu_{\hat\bbT_x} \bigl((
\mu_{\hat\bbT_x\setminus\Delta_x
}g)^2 \bigr) \bigr)
\\
&=&\delta(\ell)\mu_x \bigl(\mu_{\hat\bbT_x} \bigl((
\mu_{\hat\bbT
_x\setminus\Delta_x }g)^2 \bigr) \bigr),\nonumber
\end{eqnarray}
where $\delta(\ell):=\mu_{\hat\bbT_x}(1-\bar c_x^{(\ell)})$. Above
we used Cauchy--Schwarz to obtain the second inequality. The last
equality holds because $\mu_{\hat\bbT_x}(1-\bar c_x^{(\ell)})$ does
not depend on~$\eta_x$. Notice that $\delta(\ell)$ coincides with
$\bar p_\ell/p$ where $\bar p_\ell$ was defined at the
beginning of the proof of Theorem~\ref{teopc}.

Next we note that
%
\begin{equation}\label{last}\qquad
\mu_x \bigl(\mu_{\hat\bbT_x} \bigl(( \mu_{\hat
\bbT
_x\setminus\Delta_x }g)^2
\bigr) \bigr)=\mu_{x\cup\Delta_x
} \bigl(( \mu_{\hat\bbT_x\setminus\Delta_x }g)^2
\bigr)=\Var_{x\cup\Delta_x } ( \mu_{\hat\bbT_x\setminus\Delta_x
}g ),
\end{equation}
where we use the fact that $\mu_{x\cup\Delta_x } (\mu_{\hat
\bbT_x\setminus\Delta_x }g )=\mu_{\bbT_x}(g)=0$ by the
definition of $g$.
Then by using (\ref{A}), (\ref{least}) and (\ref{last}) we get
%
\begin{eqnarray}
\Var_x \bigl(\mu_{\hat\bbT_x} \bigl(\bigl(1-\bar
c_x^{(\ell)}\bigr)g \bigr) \bigr)&\leq&\delta(\ell)\sum
_{z\in x\cup\Delta_x}\mu_{x\cup
\Delta_x}\bigl(\Var_z(
\mu_{\hat\bbT_z} \mu_{\hat\bbT_x\setminus
\Delta_x }g)\bigr)
\nonumber\\[-8pt]\\[-8pt]
&\leq&\delta(\ell)\sum_{z\in x\cup\Delta_x}
\mu_{x\cup\Delta
_x}\bigl(\Var_z( \mu_{\hat\bbT_z}g)\bigr),\nonumber
\end{eqnarray}
where we use the convexity of the variance to obtain the second inequality.
In conclusion,
%
\begin{eqnarray}\label{F}\quad
&&
\sum_{x\in\bbT} \mu_\bbT \bigl(
\Var_x \bigl(\mu_{\hat\bbT
_x}(f) \bigr) \bigr)
\nonumber
\\
&&\qquad\le2 \sum_{x\in\bbT}\mu_\bbT \bigl(\bar
c_x^{(\ell)}\Var_x(f) \bigr) + 2\delta(\ell)\sum
_{x\in\bbT}\sum_{z\in x\cup\D_x}
\mu_\bbT \bigl(\Var_z \bigl(\mu_{\hat\bbT_z}(f) \bigr)
\bigr)
\\
&&\qquad\le4 \cD_\bbT^{(\ell)}(f) + 2(\ell+1) \delta(\ell) \sum
_{x\in
\bbT} \mu_\bbT \bigl(\Var_x
\bigl(\mu_{\hat\bbT_x}(f) \bigr) \bigr),\nonumber
\end{eqnarray}
where the factor $\ell+1$ accounts for the number of vertices $x$ such
that a
given vertex $z$ falls inside $\D_x$.\vadjust{\goodbreak}

We now appeal to Remark
\ref{bdecay} and conclude that for any $p<p_c$ there exists
$\ell_0$ (which depends
on $p$ and it diverges as $p\uparrow p_c$) such that $(\ell+1)
\delta(\ell)\le1/4$ for any $\ell\ge\ell_0$. With this choice
and recalling (\ref{A}), the Poincar\'e inequality (\ref{PI}) with
$\lambda=8$ follows uniformly in
the depth $n$ of $\bbT$. In other words the auxiliary long range model
has a positive spectral gap greater than $1/8$ if $\ell\ge\ell_0$.

We are now in a position to conclude the proof in the oriented
case. Starting from (\ref{PI}) and using path arguments exactly as in
Section 5 of~\cite{CMRT}, we conclude that, for any $\ell\ge\ell_0$ we
can find a constant $\lambda(\ell,k,j)\ge1$ independent of $n$ such that
\[
\Var_\bbT(f)\le\lambda(\ell,k,j)\sum_{x\in\bbT}
\mu_\bbT \bigl(\bar c_{\bbT,x} \Var_x(f) \bigr).
\]
Thus, thanks to Remark~\ref{finite-infinite}, we can conclude that
the spectral gap of the oriented model on the infinite tree $\bar\bbT^k$ is bounded from below by
$\lambda(\ell,k,j)^{-1}$.

\begin{remark}
The dependence on $p$ of
$\lambda(\ell,k,j)$ comes from the fact that $\ell>\ell_0(p,j,k)$. Clearly
the critical scale $\ell_0$ diverges as $p\uparrow p_c$.
\end{remark}
\noqed\end{pf*}
\begin{pf*}{Proof of Theorem~\ref{main}: the unoriented case}
For an arbitrary vertex $r\in\bbT^k$ we introduce
an auxiliary \textit{block dynamics}, reversible w.r.t. the measure
$\mu$, as
follows. With rate one the block chain resamples the current
configuration in
$\bbT^k\setminus r$ from the equilibrium measure, and, always with
rate one, it
resamples the variable $\h_r$ if and only if the constraint at the
root is satisfied [i.e., $c_r(\h)=1$].

For such auxiliary block chain it is easy to prove a Poincar\'e
inequality of
the form (compare to Proposition 4.4 in~\cite{CMRT})
%
\begin{equation}
\label{FF} \Var(f) \le\g\mu \bigl(c_r \Var_r(f) +
\Var_{\bbT^k\setminus r
}(f) \bigr)
\end{equation}
for some constant $\g=\g(j,k)\ge1$.\vspace*{1pt}

Observe now that $\bbT^k\setminus r$ is the union of $k+1$ copies of
the rooted tree $\bar
\bbT^k$ so that
\[
\Var_{\bbT^k\setminus r }(f) \le\sum_{y\in\cN_r}
\mu_{\bbT
^k\setminus r }\bigl(\Var_{\bbT_y^k}(f)\bigr).
\]
Thanks to the result in the oriented case and using $\bar c_x \le c_x$,
we get
%
\begin{equation}
\label{F1} \Var_{\bbT_y^k}(f)\le\lambda\sum_{x\in
\bbT^k_y}
\mu_{\bbT^k_y} \bigl(\bar c_x \Var_x(f) \bigr) \le
\lambda\sum_{x\in
\bbT^k_y}\mu_{\bbT^k_y}
\bigl(c_x \Var_x(f) \bigr),
\end{equation}
where $\lambda=\lambda(\ell,k,j)$.
Thus
%
\begin{equation}
\label{F2} \mu \bigl(\Var_{\bbT^k\setminus r }(f) \bigr)\le\lambda\mathop{\sum
_{x\in\bbT^k}}_{x\neq r}\mu\bigl(c_x\Var_x(f)\bigr).
\end{equation}
Inserting (\ref{F2}) into (\ref{FF}) we conclude that the spectral
gap of
the FA-jf model is bounded below by $(\g\lambda)^{-1}$.\vadjust{\goodbreak}
\end{pf*}
\begin{pf*}{Proof of Corollary~\ref{Coro}}
We closely follow the proof of a similar result given in \cite
{Martinelli-Wouts}. Recall that $\bbT$ is the finite sub-tree
consisting of the first $n$
levels of $\bar\bbT^k$ and that $h_t^\eta(\s)$ denotes the relative
density w.r.t.
$\mu_\bbT$ of the law at time $t$ of the oriented chain started at
$\eta$.
We can then write
\[
h_{t+s}^\eta(\cdot)= e^{t \cL_{\bbT}}\bigl(h_s^\eta
\bigr) (\cdot)
\]
together with
\[
h_s^\eta(\sigma)= \frac{\nu^\eta_s(\sigma_r \tc
\bigcap_{y\in\cK_r} \{\sigma_{\bbT_y}\})}{\mu_\bbT(\sigma_r)}\prod
_{y
\in\cK_r}h_s^\eta(\sigma_{\bbT_y})
\le\frac{1}{\min(p,q)}\prod_{y \in\cK_r}h_s^\eta(
\sigma_{\bbT_y}).
\]
Above $\bbT_y$ denotes the sub-tree of $\bbT$ rooted at $y$ and we
used the fact that, because of the orientation of the model, the
marginal of $\nu_t^\eta$ on $\O_{\bbT\setminus
r}$ is the product over $y\in\cK_r$ of its marginals on $\O_{\bbT_y}$.
Therefore
%
\begin{eqnarray}\label{pippo}
\Var_\bbT \bigl(h^\eta_{t+s} \bigr)&=&
\Var_\bbT \bigl( e^{t\bar\cL
_{\bbT}} h_s^\eta
\bigr) \le e^{-\operatorname{gap}(\bar\cL_\bbT) t} \Var_\bbT \bigl(h_s^\eta
\bigr)
\nonumber
\\
&\le& e^{-\operatorname{gap}(\bar\cL_\bbT) t} \frac{1}{\min(p,q)^{2}} \prod_{y \in\cK_r}
\mu_{\bbT_y} \bigl(\bigl[ h_s^\eta
\bigr]^2 \bigr)
\\
&=& e^{-\operatorname{gap}(\bar\cL_\bbT) t} \frac{1}{\min(p,q)^{2}} \prod_{y
\in\cK_r}
\bigl(\Var_{\bbT_y}\bigl(h_s^\eta\bigr)+1 \bigr).\nonumber
\end{eqnarray}
Let now $t_n:=\inf\{t\ge0\dvtx  \max_\eta\Var_\bbT^\eta(h_t^\eta
)\le
1/4\}$ so that, by definition, $t_n=T_2$. If in (\ref{pippo}) we
choose $s=t_{n-1}$ we get
\[
\Var_\bbT \bigl(h^\eta_{t+t_{n-1}} \bigr)\le
e^{-\operatorname{gap}(\bar\cL
_\bbT)
t}\frac{5^k}{4^k{\min(p,q)}^2},
\]
because each sub-tree $\bbT_y$ has $n-1$ levels.

Thus, if $t^*$ is so large that $\frac{5^k}{4^k{\min(p,q)}^2}
e^{-\operatorname{gap}(\bar\cL_\bbT) t^*} \le1/4$,
then
\[
\max_\eta\Var_\bbT\bigl(h^\eta_{t^*+t_{n-1}}
\bigr)\le1/4,
\]
that is, $T_2=t_n\le t^*+t_{n-1}\le\cdots\leq t^* n$. That completes the
proof of the upper bound.

The linear lower bound, $T_1\ge c n$ for some constant $c >0$, follows
immediately from the fact that,
starting from the configuration $\eta$ with $\eta_x=1 \ \forall x\in
\bbT$, routine bounds show that the influence from the leaves cannot
propagate faster than
linear in time; see, for example,~\cite{Martinelli97}.
\end{pf*}

\subsection{The super-critical phase $p>p_c$}

\mbox{}
\begin{pf*}{Proof of Proposition~\ref{soprapc}}
If $p>p_c$, then with positive probability the root $r$ belongs to an infinite
cluster of occupied vertices which is stable upon iterations of the\vadjust{\goodbreak}
bootstrap map $\bar B$. Clearly any vertex belonging
to such a cluster can never change its occupation variable during the
dynamics of OFA-jf. Hence the result. The result for the nonoriented
model can be established via the same lines by replacing the root $r$
with any arbitrary vertex of $\bbT^k$.
\end{pf*}
\begin{pf*}{Proof of Theorem~\ref{abovepc}}
Fix $n$ and consider for simplicity only the rooted case, the unrooted
one being treated along the same lines. As before we denote by $\bbT$ the
$k$-ary rooted tree of depth $n$ and root $r$. We begin by proving the
stated upper bound.

Choose as test function $f$ to be used in the Poincar\'e
inequality
\[
\gap(\cL_\bbT)\le\cD_\bbT(f)/\Var_\bbT(f)
\qquad\forall f\in L^2(\bbT)
\]
the
indicator of the event $A$ that the root is occupied after
$n-1$ iterations of the bootstrap map $\bar B$. If $p>p_c$, then $\Var_\bbT(f)>0$
uniformly in $n$.

Next we compute the Dirichlet form $\cD_\bbT(f)$. We first observe
that if $x\in\bbT$ is not a leaf of $\bbT$, then
the corresponding contribution $\mu_\bbT(\bar c_{\bbT,x} \Var_x(f))$ to the Dirichlet
form vanishes. Otherwise one could connect $A$ to $\O_\bbT\setminus A$
by means of
a legal flip, that is, one with $\bar c_{\bbT,x}=1$. But that is clearly
impossible by the definition of $A$. If instead $x$ is a leaf of
$\bbT$, so that $\bar c_{\bbT,x}\equiv1$ by definition, then
\[
\mu_\bbT\bigl(\Var_x(f)\bigr)= 2\mu_\bbT
\bigl(\eta\in A; \eta^x\notin A\bigr).
\]
The latter probability can be computed explicitly and it is equal to
$\prod_{y\preceq x}p_y$
where $y\preceq x$ means that $y$ is an ancestor of $x$, and $p_y$ is
the probability that $y$ is occupied and that exactly $j-1$ out of the
$k-1$ children of $y$ which are not ancestors of $x$ are not occupied after
$n-d_y-1$ iteration of the bootstrap map $\bar B$. Since the probability
$p^{(n)}$ that the root is occupied after $n$-iterations of the
bootstrap map converges exponentially fast to the largest fixed point
$p_{\infty}$ of the map $g_p(\cdot)$ defined in Theorem~\ref{main},
we get that
\[
\prod_{y\preceq x}p_y \le C \biggl(p \pmatrix{k-1
\cr j-1} (1-p_{\infty
})^{j-1}p_\infty^{k-j}
\biggr)^n
\]
for some positive constant $C$. In conclusion
\[
\cD_\bbT(f)\le C \biggl(kp \pmatrix{k-1\cr j-1} p_{\infty}^{k-j}(1-p_\infty)^{j-1}
\biggr)^n.
\]
The proof of the upper bound is complete once we observe that
\[
kp \pmatrix{k-1\cr j-1} p_{\infty}^{k-j}(1-p_\infty)^{j-1}=
\frac
{d}{d\lambda}g_p(\lambda) \bigg|_{\lambda=p_\infty}<1.
\]
We now turn to the lower bound. The proof is based on the same
argument used in Theorem~\ref{teogap} to treat the unoriented case
which we now shortly detail.

By monotonicity of the rates as functions of $j$ we may assume $j=k$.
As before,
consider
the auxiliary block dynamics in which:
\begin{itemize}
\item each sub-tree rooted at one of
the children of the root with rate one updates at the same time all its
vertices by choosing the new configuration from the equilibrium
distribution;
\item the root $r$, with rate one and if and only if all its children are
empty, refreshes its occupation variable
by sampling a new value from the equilibrium measure.
\end{itemize}
It is easy to check
that, for any $p\in(0,1)$, the spectral gap of the block dynamics is
positive uniformly in $n$ so that a uniform Poincar\'e inequality holds
\[
\Var_\bbT(f)\le C \mu_\bbT \biggl(c_r
\Var_r(f) + \sum_{x\in
\cN_r}
\Var_{\bar\bbT_x}(f) \biggr) \qquad\forall f
\]
for some $C>0$ independent of $n$.

For notational convenience let $\g(n):=\gap(\bar\cL_\bbT)^{-1}$.
By definition, for each \mbox{$x\in\cN_r$}, $\Var_{\bar\bbT_x}(f)\le
\g(n-1)\cD_{\bar\bbT_x}(f)$. Therefore
\[
\Var_\bbT(f)\le C\max\bigl(1,\g(n-1)\bigr)\cD_\bbT(f),
\]
that is,
\[
\g(n)\le C\max\bigl(1,\g(n-1)\bigr)\le\cdots\le C^n.
\]
\upqed\end{pf*}

\section{Extensions to KCSM on $\bbZ^d$}
\label{generalizations}

In this section we discuss some applications of the technique that we
have devised to prove Theorem~\ref{teogap}.
We show in particular that this technique allows us to recover the
positivity of the spectral gap in the whole ergodicity region for the
KCSM on $\bbZ^d$ which were studied in~\cite{CMRT} via completely
different methods.
We start by treating explicitly the case of the north--east model, and
then we will describe how to extend the analysis to more general models
\cite{RJM}.
%
\begin{definition}
The North--East (N--E) model is a KCSM on $\bbZ^2$ for
which the constraint at $x\in\bbZ^2$ requires the northern \textit{and}
eastern neighbor of $x$ to be empty. More precisely it is a continuous
time Markov process on $\O=\{0,1\}^{\bbZ^2}$ with generator $\cL$
defined as in
Definition~\ref{defFA} but with the sum in the generator now running
on the sites of $\bbZ^2$ and with constraints
%
\begin{equation}
\label{cxNE} c_{x}(\eta)= %
\cases{1, &\quad if $
\eta_{x+\vec e_1}= \eta_{x+\vec e_2}=0$,
\cr
0, &\quad otherwise, } %
\end{equation}
with $\vec e_1$ and $\vec e_2$ the Euclidean unit vectors on $\bbZ^2$.
\end{definition}

Let us recall some well-known properties of the North--East model
\cite{KL,CMRT} (in particular we refer the reader to Section 6.4 of
\cite{CMRT}
where these results have been derived by using the analog of our
Proposition~\ref{propo} and via the results on oriented percolation of
\cite{Schonmann} and~\cite{Durrett}).

Let $p_c$ be the critical density defined as in (\ref{eq1}), and let
the associated bootstrap map $B$ be defined exactly as in
(\ref{eqbootstrapmap}). Let $\mu^{(n)}$ be the measure obtained by
iterating $B$ $n$-times starting from $\mu$, and call $p_n$ be the
probability that a vertex is occupied under $\mu^{(n)}$.
%
\begin{MainProposition}\label{propNE}
$p_c$ coincides with the critical threshold for
oriented percolation in
$\bbZ^2$. In particular~\cite{Durrett} $p_c\in(0,1)$.
Moreover, for any $p<p_c$,
\[
\lim_{n\to\infty}n^2 p_{n}=0.
\]
\end{MainProposition}
We will now prove via the technique described in Section~\ref{subcrit}
the following result.
%
\begin{MainTheorem}\label{gapNE}
Assume $p<p_c$. Then $\gap(\cL)>0$.
\end{MainTheorem}
\begin{pf}
As for the models on trees, we prove a lower bound on the spectral
gap on an arbitrarily large finite region $\L$ of $\bbZ^2$ with
proper boundary conditions
which is uniform in the size of the region. Then the result on
infinite volume follows by a standard limiting procedure.

The finite region $\L\subset\bbZ^2$ that we consider consists of all
the points of
$\bbZ^2$ inside the right triangle
$\Lambda\subset\bbR^2$ with a vertex in the origin, a vertex in
$n\vec e_1$ and a vertex in $n\vec e_2$ where $n$ is a large
integer. We will consider the North--East model in $\L$ with empty
boundary conditions, namely with $\L$-dependent constraints $c_{\L,x}$
given by $c_{\L,x}(\eta)=c_x(\eta^0)$ for any $\eta\in\O_\L$,
where $\eta^0$ is as in
Definition~\ref{finite} with the obvious modifications.

For any $x\in\Lambda$ let
\[
C_x:=\{z\in\Lambda\dvtx  z\cdot\vec e_1\geq x\cdot\vec
e_1 \mbox{ and } z\cdot\vec e_2\geq x\cdot\vec
e_2\}
\]
be the right
cone with vertex at $x$, and let $\hat C_x:=C_x\setminus x$.

Our first claim is the analog of inequality (\ref{A})
proved in the tree case. More precisely,
%
\begin{claim}
%
\begin{equation}
\label{AB} \Var_{\Lambda}(f)\leq\sum_x
\mu_{\Lambda} \bigl(\Var_x \bigl(\mu_{\hat
C_x}(f) \bigr)
\bigr) \qquad\forall f.
\end{equation}
\end{claim}
\begin{pf}
For $j=0,1,\ldots,n$ let $\L_j$ be the set of vertices in $\L$ with
$\ell_1$ distance from the origin at least $L-j$.
Then
\begin{eqnarray*}
\Var_\L(f)&=&\mu_\L \bigl(\Var_{\L_0}(f)
\bigr) +\Var_\L \bigl(\mu_{\L_0}(f) \bigr)
\\
&=& \mu_\L \bigl(\Var_{\L_0}(f) \bigr) +\mu_\L
\bigl(\Var_{\L_1} \bigl[\mu_{\L_0}(f) \bigr] \bigr) +
\Var_\L \bigl(\mu_{\L_1} \bigl[\mu_{\L_0}(f) \bigr]
\bigr)
\\
&\vdots&
\\
&=&\mu_\L \bigl(\Var_{\L_0}[f] \bigr)+\sum
_{j=0}^{n-1}\mu_\L \bigl(
\Var_{\L_{j+1}} \bigl[\mu_{\L_{j}}(f) \bigr] \bigr).
\end{eqnarray*}
Thanks to (\ref{var}),
\begin{eqnarray*}
\Var_{\L_{j+1}} \bigl[\mu_{\L_{j}}(f) \bigr]&=& \Var_{\L
_{j+1}\setminus
\L_j}
\bigl[\mu_{\L_{j}}(f) \bigr]\le \sum_{x\in\L_{j+1}\setminus\L_n}
\mu_{\L_{j+1}\setminus
\L_j} \bigl(\Var_x\bigl(\mu_{\L_j}(f) \bigr)\bigr)
\\
&\le&\sum_{x\in\L_{j+1}\setminus\L_j}\mu_{\L_{j+1}\setminus
\L_j} \bigl(
\Var_x\bigl(\mu_{\hat C_x}(f) \bigr)\bigr),
\end{eqnarray*}
where in the last inequality we used the fact that $\hat C_x\subset\L_j$
for all $x\in\L_{j+1}$ together with the standard convexity property of
the variance. Analogously,
\[
\mu_\L \bigl(\Var_{\L_0}[f] \bigr)\le\sum
_{x\in\L_0}\mu_\L \bigl(\Var_x(f) \bigr).
\]
The proof of the claim is complete if we observe that, for any $f$,
\[
\mu_\L \bigl(\mu_{\L_{j+1}\setminus
\L_j}(f) \bigr)=\mu_\L(f).
\]
\upqed\end{pf}

Back to the proof of Theorem~\ref{gapNE}, for any integer $\ell\ll
n$ let $c_x^{\ell}(\eta)$ be defined exactly as the long-range
constraints $\bar c_x^{(\ell)}$ given in Definition~\ref{meso} with
the tree $\bbT$ replaced by the region $\L$ and $\hat\bbT_x$
replaced by
$\hat C_x$.
Then, by using the key inequality (\ref{AB}) and by following exactly
the same route of the proof of Theorem~\ref{teogap}, we obtain
%
\begin{eqnarray}
\label{FNE}\hspace*{28pt} \Var_\L(f)&\le&\sum_x
\mu_{\Lambda} \bigl(\Var_x \bigl(\mu_{\hat
C_x}(f) \bigr)
\bigr)
\\
&\le&4 \sum_{x\in\L}\mu_\L
\bigl(c_x^{(\ell)}\Var_x(f) \bigr) + \frac2p(
\ell+1)^2 p_{\ell}\sum_x
\mu_\L \bigl(\Var_x \bigl(\mu_{\hat\bbT_x}(f) \bigr)
\bigr),
\end{eqnarray}
where the factor $(\ell+1)^2$ [instead of $(\ell+1)$ of (\ref{F})]
accounts for the number of vertices $x$ such
that their $\ell_1$-distance from a given vertex $z$ is at most $\ell$.
Proposition~\ref{propNE} implies that there exists $\ell_0=\ell_0(p)$ such that
\[
\frac2p (\ell+1)^2p_\ell<1/2 \qquad\forall\ell\geq
\ell_0.
\]
Therefore, if $\ell\ge\ell_0$,
\[
\Var_\L(f) \le8 \sum_{x\in
\L}
\mu_\L \bigl(c_x^{(\ell)}\Var_x(f)
\bigr).
\]
Elementary path arguments (see also~\cite{CMRT}) show now that
\[
\sum_{x\in
\L}\mu_\L \bigl(c_x^{(\ell)}
\Var_x(f) \bigr) \le C(\ell) \sum_{x\in
\L}
\mu_\L \bigl(c_{\L,x}\Var_x(f) \bigr)
\]
for some finite constant $C$ independent of $n$. The proof is complete.
\end{pf}
%
\begin{remark}
Via a proper generalization of our technique we can establish the
positivity of the spectral gap for all the KCSM\vadjust{\goodbreak} covered by Theorem~3.3
of~\cite{CMRT}. These include, besides N--E model, some of the
KCSM which have been most studied in physics literature, namely the
East model on $\bbZ$, the Friedrickson--Andersen model on $\bbZ^d$ and
the modified basic model on $\bbZ^d$; see Section 2.3 of
\cite{CMRT} for the definitions. More precisely our technique allows us
to prove Theorem 4.1 of~\cite{CMRT} (in a completely different way),
namely to establish the positivity of the spectral gap in a proper
regime for the so-called \textit{*-general model}
\cite{CMRT}. Then the proof of positivity of the spectral gap for each
specific KCSM can be completed via the renormalization technique
detailed in Section 5 of~\cite{CMRT}.
Along the same lines we can also recover the positivity of the spectral
gap for the spiral model, a result which was previously established in
\cite{LNM}.
\end{remark}

\section*{Acknowledgments}
We thank the Laboratoire de Probabilit\'{e}s et Mod\`{e}les
Al\'{e}atoires, the University Paris VII and the Department of
Mathematics of the University of Roma Tre for the support and the kind
hospitality.


%
%

\printaddresses


\begin{thebibliography}{26}

\bibitem{Aldous}
%
\begin{barticle}[mr]
\bauthor{\bsnm{Aldous},~\bfnm{David}\binits{D.}} \AND
\bauthor{\bsnm{Diaconis},~\bfnm{Persi}\binits{P.}}
(\byear{2002}).
\btitle{The asymmetric one-dimensional constrained {I}sing model: Rigorous
results}.
\bjournal{J. Stat. Phys.}
\bvolume{107}
\bpages{945--975}.
\bid{doi={10.1023/A:1015170205728}, issn={0022-4715}, mr={1901508}}
\bptok{imsref}%
\end{barticle}
%
\endbibitem

\bibitem{Peres}
%
\begin{barticle}[mr]
\bauthor{\bsnm{Balogh},~\bfnm{J{\'o}zsef}\binits{J.}},
\bauthor{\bsnm{Peres},~\bfnm{Yuval}\binits{Y.}} \AND
\bauthor{\bsnm{Pete},~\bfnm{G{\'a}bor}\binits{G.}}
(\byear{2006}).
\btitle{Bootstrap percolation on infinite trees and non-amenable groups}.
\bjournal{Combin. Probab. Comput.}
\bvolume{15}
\bpages{715--730}.
\bid{doi={10.1017/S0963548306007619}, issn={0963-5483}, mr={2248323}}
\bptok{imsref}%
\end{barticle}
%
\endbibitem

\bibitem{CMRT}
%
\begin{barticle}[mr]
\bauthor{\bsnm{Cancrini},~\bfnm{N.}\binits{N.}},
\bauthor{\bsnm{Martinelli},~\bfnm{F.}\binits{F.}},
\bauthor{\bsnm{Roberto},~\bfnm{C.}\binits{C.}} \AND
\bauthor{\bsnm{Toninelli},~\bfnm{C.}\binits{C.}}
(\byear{2008}).
\btitle{Kinetically constrained spin models}.
\bjournal{Probab. Theory Related Fields}
\bvolume{140}
\bpages{459--504}.
\bid{doi={10.1007/s00440-007-0072-3}, issn={0178-8051}, mr={2365481}}
\bptok{imsref}%
\end{barticle}
%
\endbibitem

\bibitem{LNM}
%
\begin{bincollection}[mr]
\bauthor{\bsnm{Cancrini},~\bfnm{N.}\binits{N.}},
\bauthor{\bsnm{Martinelli},~\bfnm{F.}\binits{F.}},
\bauthor{\bsnm{Roberto},~\bfnm{C.}\binits{C.}} \AND
\bauthor{\bsnm{Toninelli},~\bfnm{C.}\binits{C.}}
(\byear{2009}).
\btitle{Facilitated spin models: Recent and new results}.
In \bbooktitle{Methods of Contemporary Mathematical Statistical Physics}
(\beditor{R. Kotecky}, ed.).
\bseries{Lecture Notes in Math.}
\bvolume{1970}
\bpages{307--340}.
\bpublisher{Springer}, \blocation{Berlin}.
\bid{mr={2581609}}
\bptnote{check year}%
\bptok{imsref}%
\end{bincollection}
%
\endbibitem

\bibitem{critical-tree}
%
\begin{bmisc}[auto:STB|2013/01/29|08:09:18]
\bauthor{\bsnm{Cancrini},~\bfnm{N.}\binits{N.}},
\bauthor{\bsnm{Martinelli},~\bfnm{F.}\binits{F.}},
\bauthor{\bsnm{Roberto},~\bfnm{C.}\binits{C.}} \AND
\bauthor{\bsnm{Toninelli},~\bfnm{C.}\binits{C.}}
(\byear{2012}).
\bhowpublished{Mixing time of a kinetically constrained spin model on trees: Power law
scaling at criticality. Preprint}.
\bptok{imsref}%
\end{bmisc}
%
\endbibitem

\bibitem{CMST}
%
\begin{barticle}[mr]
\bauthor{\bsnm{Cancrini},~\bfnm{N.}\binits{N.}},
\bauthor{\bsnm{Martinelli},~\bfnm{F.}\binits{F.}},
\bauthor{\bsnm{Schonmann},~\bfnm{R.}\binits{R.}} \AND
\bauthor{\bsnm{Toninelli},~\bfnm{C.}\binits{C.}}
(\byear{2010}).
\btitle{Facilitated oriented spin models: Some non equilibrium results}.
\bjournal{J. Stat. Phys.}
\bvolume{138}
\bpages{1109--1123}.
\bid{doi={10.1007/s10955-010-9923-x}, issn={0022-4715}, mr={2601425}}
\bptok{imsref}%
\end{barticle}
%
\endbibitem

\bibitem{chalupa}
%
\begin{barticle}[auto:STB|2013/01/29|08:09:18]
\bauthor{\bsnm{Chalupa},~\bfnm{J.}\binits{J.}},
\bauthor{\bsnm{Leath},~\bfnm{P.~L.}\binits{P.~L.}} \AND
\bauthor{\bsnm{Reich},~\bfnm{G.~R.}\binits{G.~R.}}
(\byear{1979}).
\btitle{Bootstrap percolation on a Bethe lattice}.
\bjournal{J. Phys. C: Solid State Phys.}
\bvolume{12}
\bpages{L31--L35}.
\bptok{imsref}%
\end{barticle}
%
\endbibitem

\bibitem{Durrett}
%
\begin{barticle}[mr]
\bauthor{\bsnm{Durrett},~\bfnm{Richard}\binits{R.}}
(\byear{1984}).
\btitle{Oriented percolation in two dimensions}.
\bjournal{Ann. Probab.}
\bvolume{12}
\bpages{999--1040}.
\bid{issn={0091-1798}, mr={0757768}}
\bptok{imsref}%
\end{barticle}
%
\endbibitem

\bibitem{FS}
%
\begin{barticle}[mr]
\bauthor{\bsnm{Fontes},~\bfnm{Luiz~Renato}\binits{L.~R.}} \AND
\bauthor{\bsnm{Schonmann},~\bfnm{Roberto~H.}\binits{R.~H.}}
(\byear{2008}).
\btitle{Threshold {$\theta\geq2$} contact processes on homogeneous trees}.
\bjournal{Probab. Theory Related Fields}
\bvolume{141}
\bpages{513--541}.
\bid{doi={10.1007/s00440-007-0092-z}, issn={0178-8051}, mr={2391163}}
\bptok{imsref}%
\end{barticle}
%
\endbibitem

\bibitem{Fredrickson1}
%
\begin{barticle}[auto:STB|2013/01/29|08:09:18]
\bauthor{\bsnm{Fredrickson},~\bfnm{G.}\binits{G.}} \AND
\bauthor{\bsnm{Andersen},~\bfnm{H.}\binits{H.}}
(\byear{1984}).
\btitle{Kinetic Ising model of the Glass transition}.
\bjournal{Phys. Rev. Lett.}
\bvolume{53}
\bpages{1244--1247}.
\bptok{imsref}%
\end{barticle}
%
\endbibitem

\bibitem{Fredrickson2}
%
\begin{barticle}[auto:STB|2013/01/29|08:09:18]
\bauthor{\bsnm{Fredrickson},~\bfnm{G.}\binits{G.}} \AND
\bauthor{\bsnm{Andersen},~\bfnm{H.}\binits{H.}}
(\byear{1985}).
\btitle{Facilitated kinetic Ising models and the glass transition}.
\bjournal{J. Chem. Phys.}
\bvolume{83}
\bpages{5822--5831}.
\bptok{imsref}%
\end{barticle}
%
\endbibitem

\bibitem{GST}
%
\begin{bbook}[auto:STB|2013/01/29|08:09:18]
\bauthor{\bsnm{Garrahan},~\bfnm{J.}\binits{J.}},
\bauthor{\bsnm{Sollich},~\bfnm{P.}\binits{P.}} \AND
\bauthor{\bsnm{Toninelli},~\bfnm{C.}\binits{C.}}
(\byear{2011}).
\btitle{Dynamical Heterogeneities in Glasses, Colloids, and Granular Media}.
\bpublisher{Oxford Univ. Press}, \blocation{Oxford}.
\bnote{Available at arXiv:\arxivurl{1009.6113}}.
\bptok{imsref}%
\end{bbook}
%
\endbibitem

\bibitem{Grimmett}
%
\begin{bbook}[mr]
\bauthor{\bsnm{Grimmett},~\bfnm{Geoffrey}\binits{G.}}
(\byear{1999}).
\btitle{Percolation},
\bedition{2nd} ed.
\bseries{Grundlehren der Mathematischen Wissenschaften [Fundamental Principles
of Mathematical Sciences]}
\bvolume{321}.
\bpublisher{Springer}, \blocation{Berlin}.
\bid{mr={1707339}}
\bptok{imsref}%
\end{bbook}
%
\endbibitem

\bibitem{RJM}
%
\begin{barticle}[auto:STB|2013/01/29|08:09:18]
\bauthor{\bsnm{Jackle},~\bfnm{J.}\binits{J.}},
\bauthor{\bsnm{Mauch},~\bfnm{F.}\binits{F.}} \AND
\bauthor{\bsnm{Reiter},~\bfnm{J.}\binits{J.}}
(\byear{1992}).
\btitle{Blocking transitions in lattice spin models with directed kinetic
constraints}.
\bjournal{Phys. A}
\bvolume{184}
\bpages{458--476}.
\bptok{imsref}%
\end{barticle}
%
\endbibitem

\bibitem{KL}
%
\begin{barticle}[mr]
\bauthor{\bsnm{Kordzakhia},~\bfnm{George}\binits{G.}} \AND
\bauthor{\bsnm{Lalley},~\bfnm{Steven~P.}\binits{S.~P.}}
(\byear{2006}).
\btitle{Ergodicity and mixing properties of the northeast model}.
\bjournal{J. Appl. Probab.}
\bvolume{43}
\bpages{782--792}.
\bid{doi={10.1239/jap/1158784946}, issn={0021-9002}, mr={2274800}}
\bptok{imsref}%
\end{barticle}
%
\endbibitem

\bibitem{Liggett}
%
\begin{bbook}[mr]
\bauthor{\bsnm{Liggett},~\bfnm{Thomas~M.}\binits{T.~M.}}
(\byear{1985}).
\btitle{Interacting Particle Systems}.
\bseries{Grundlehren der Mathematischen Wissenschaften [Fundamental Principles
of Mathematical Sciences]}
\bvolume{276}.
\bpublisher{Springer}, \blocation{New York}.
\bid{doi={10.1007/978-1-4613-8542-4}, mr={0776231}}
\bptok{imsref}%
\end{bbook}
%
\endbibitem

\bibitem{Martinelli97}
%
\begin{bincollection}[mr]
\bauthor{\bsnm{Martinelli},~\bfnm{Fabio}\binits{F.}}
(\byear{1999}).
\btitle{Lectures on {G}lauber dynamics for discrete spin models}.
In \bbooktitle{Lectures on Probability Theory and Statistics ({S}aint-{F}lour,
1997)}.
\bseries{Lecture Notes in Math.}
\bvolume{1717}
\bpages{93--191}.
\bpublisher{Springer}, \blocation{Berlin}.
\bid{doi={10.1007/978-3-540-48115-7_2}, mr={1746301}}
\bptok{imsref}%
\end{bincollection}
%
\endbibitem

\bibitem{MSW}
%
\begin{barticle}[mr]
\bauthor{\bsnm{Martinelli},~\bfnm{Fabio}\binits{F.}},
\bauthor{\bsnm{Sinclair},~\bfnm{Alistair}\binits{A.}} \AND
\bauthor{\bsnm{Weitz},~\bfnm{Dror}\binits{D.}}
(\byear{2004}).
\btitle{Glauber dynamics on trees: Boundary conditions and mixing time}.
\bjournal{Comm. Math. Phys.}
\bvolume{250}
\bpages{301--334}.
\bid{doi={10.1007/s00220-004-1147-y}, issn={0010-3616}, mr={2094519}}
\bptok{imsref}%
\end{barticle}
%
\endbibitem

\bibitem{Martinelli-Wouts}
%
\begin{barticle}[mr]
\bauthor{\bsnm{Martinelli},~\bfnm{Fabio}\binits{F.}} \AND
\bauthor{\bsnm{Wouts},~\bfnm{Marc}\binits{M.}}
(\byear{2012}).
\btitle{Glauber dynamics for the quantum {I}sing model in a transverse
field on
a regular tree}.
\bjournal{J. Stat. Phys.}
\bvolume{146}
\bpages{1059--1088}.
\bid{doi={10.1007/s10955-012-0436-7}, issn={0022-4715}, mr={2902454}}
\bptok{imsref}%
\end{barticle}
%
\endbibitem

\bibitem{Ritort-Sollich}
%
\begin{barticle}[auto:STB|2013/01/29|08:09:18]
\bauthor{\bsnm{Ritort},~\bfnm{F.}\binits{F.}} \AND
\bauthor{\bsnm{Sollich},~\bfnm{P.}\binits{P.}}
(\byear{2003}).
\btitle{Glassy dynamics of kinetically constrained models}.
\bjournal{Adv. Phys.}
\bvolume{52}
\bpages{219--342}.
\bptok{imsref}%
\end{barticle}
%
\endbibitem

\bibitem{Gine}
%
\begin{bincollection}[mr]
\bauthor{\bsnm{Saloff-Coste},~\bfnm{Laurent}\binits{L.}}
(\byear{1997}).
\btitle{Lectures on finite {M}arkov chains}.
In \bbooktitle{Lectures on Probability Theory and Statistics ({S}aint-{F}lour,
1996)}.
\bseries{Lecture Notes in Math.}
\bvolume{1665}
\bpages{301--413}.
\bpublisher{Springer}, \blocation{Berlin}.
\bid{doi={10.1007/BFb0092621}, mr={1490046}}
\bptok{imsref}%
\end{bincollection}
%
\endbibitem

\bibitem{STBT}
%
\begin{barticle}[mr]
\bauthor{\bsnm{Sausset},~\bfnm{Fran{\c{c}}ois}\binits{F.}},
\bauthor{\bsnm{Toninelli},~\bfnm{Cristina}\binits{C.}},
\bauthor{\bsnm{Biroli},~\bfnm{Giulio}\binits{G.}} \AND
\bauthor{\bsnm{Tarjus},~\bfnm{Gilles}\binits{G.}}
(\byear{2010}).
\btitle{Bootstrap percolation and kinetically constrained models on hyperbolic
lattices}.
\bjournal{J. Stat. Phys.}
\bvolume{138}
\bpages{411--430}.
\bid{doi={10.1007/s10955-009-9903-1}, issn={0022-4715}, mr={2594903}}
\bptok{imsref}%
\end{barticle}
%
\endbibitem

\bibitem{Schonmann}
%
\begin{barticle}[mr]
\bauthor{\bsnm{Schonmann},~\bfnm{Roberto~H.}\binits{R.~H.}}
(\byear{1992}).
\btitle{On the behavior of some cellular automata related to bootstrap
percolation}.
\bjournal{Ann. Probab.}
\bvolume{20}
\bpages{174--193}.
\bid{issn={0091-1798}, mr={1143417}}
\bptok{imsref}%
\end{barticle}
%
\endbibitem

\bibitem{schwartz}
%
\begin{barticle}[auto:STB|2013/01/29|08:09:18]
\bauthor{\bsnm{Schwartz},~\bfnm{J.~M.}\binits{J.~M.}},
\bauthor{\bsnm{Liu},~\bfnm{A.~J.}\binits{A.~J.}} \AND
\bauthor{\bsnm{Chayes},~\bfnm{L.~Q.}\binits{L.~Q.}}
(\byear{2006}).
\btitle{The onset of jamming as the sudden emergence of an infinite $k$-core
cluster}.
\bjournal{Europhysics Lett.}
\bvolume{73}
\bpages{560--566}.
\bptok{imsref}%
\end{barticle}
%
\endbibitem

\bibitem{sellitto}
%
\begin{barticle}[auto:STB|2013/01/29|08:09:18]
\bauthor{\bsnm{Sellitto},~\bfnm{M.}\binits{M.}},
\bauthor{\bsnm{Biroli},~\bfnm{G.}\binits{G.}} \AND
\bauthor{\bsnm{Toninelli},~\bfnm{C.}\binits{C.}}
(\byear{2005}).
\btitle{Facilitated spin models on Bethe lattice: Bootstrap
percolation, mode
coupling transition and glassy dynamics}.
\bjournal{Europhysics Lett.}
\bvolume{69}
\bpages{496--512}.
\bptok{imsref}%
\end{barticle}
%
\endbibitem

\end{thebibliography}
\end{document}